\documentclass[12pt,leqno]{article}
\usepackage{amsthm}
\usepackage{amsfonts}
\usepackage{amssymb}
\usepackage{amsmath}

\usepackage{epsfig}
\usepackage{graphics}

\usepackage{hyperref}

\newcommand{\R}{ { \mathbb{R} } }
\newcommand{\C}{ { \mathbb{C} } }
\newcommand{\T}{ { \mathbb{T}  } }

\newcommand{\Z}{ { \mathbb{Z} } }
\newcommand{\nz}{ { \setminus \{ 0 \} } }

\newcommand{\ph}{\varphi}
\newcommand{\vp}{\vec{\varphi}}
\newcommand{\ep}{ \varepsilon }
\newcommand{\rep}{ \sqrt{\varepsilon} }

\def\curl{\operatorname{curl}}

\newcommand{\dd}{{\rm d}}

\newcommand{\dx}{ \partial_x }
\newcommand{\dy}{ \partial_y }
\newcommand{\dz}{ \partial_z }
\newcommand{\dX}{ \partial_X }
\newcommand{\dY}{ \partial_Y }
\newcommand{\dt}{ \partial_t }
\newcommand{\dT}{ \partial_T }
\newcommand{\dth}{ \partial_\theta }

\newcommand{\tendlorsque}[2]{\mathop{\longrightarrow}\limits_{#1\rightarrow#2}}
\newcommand{\Lim}[2]{\displaystyle \lim_{#1\to#2}}

\newcommand{\ie}{\emph{i.e.} }

\newcommand{\via}{\emph{via} }

\newtheorem{theo}{Theorem}[section]

\newtheorem{defn}{Definition}[section]

\newtheorem{ex}{Example}[section]

\title{About nonlinear geometric optics} 
\author{E. Dumas}
\date{}

\begin{document}

\maketitle

\begin{center}
\small
Institut Fourier, UMR 5582 (CNRS-UJF) \\
100 rue des Math\'ematiques \\
Domaine Universitaire \\
BP 74, 38402 Saint Martin d'H\`eres - France \\
email: edumas@ujf-grenoble.fr
\end{center}

\begin{abstract}
We give an idea of the evolution of mathematical nonlinear geometric
optics from its foundation by Lax in 1957, and present applications in
various fields of mathematics and physics. 
\end{abstract}


\section{Introduction} \label{secintro}

Geometric optics goes back at least to the XVIIth Century, with Fermat, 
Snell and Descartes, who described the ``paths" (rays) followed by the light. 
Nowadays, Physics tells us that we may reasonably replace the waves from 
Quantum Mechanics with classical particles, in the semi-classical 
approximation (when considering Planck's constant $\hbar$, or the wavelength, 
as infinitely small). The mathematical transcription of these problems 
consists in studying the asymptotic behavior of solutions to partial
differential equations where different scales (represented by small
parameters) are present, often in a high frequency oscillatory context.

We present the first historical
results of the field, and then review some extensions and applications
of the method. We shall see how geometric optics applies to Maxwell's
equations (from optics, ferromagnetism, \dots), to the
wave or Klein-Gordon equation, to fluid dynamics and plasma physics, 
to general hyperbolic systems and conservation laws, as well as to
nonlinear Schr\"odinger equations, among others. Furthermore, we shall see that
it does not apply to oscillatory problems only, but also to boundary
layers, shocks and long waves problems.
\tableofcontents

Quoting some reviews and introductory texts may be useful. A review on 
nonlinear geometric optics in 1998 is due to Joly, M\'etivier and Rauch 
\cite{JMR99}, and the online book of Rauch \cite{Rau} is a nice
introduction to  
the subject. Majda \cite{Maj84} made major contributions in the 80ies. 
Good Physics textbooks on nonlinear optics are Newell and Moloney \cite{NM92} 
and Boyd \cite{Boy92}. Kalyakin \cite{Kal89} and Hunter \cite{Hun95} review 
many questions, methods and applications about these multi-scale
problems, and Whitham  \cite{Whi74} had pioneering contributions in
the modeling of nonlinear waves. 

\section{First steps} \label{secfirst}

\subsection{Linear geometric optics} 

The first rigorous result in mathematical 
geometric optics is due to Lax \cite{Lax57}, who shows that strictly 
hyperbolic systems admit \emph{WKB} (for Wentzel, Kramers and Brillouin) 
solutions. These have the \emph{phase-amplitude} representation
$$u^\ep(t,x) = a^\ep(t,x)~e^{i\ph(t,x)/\ep}.$$
The wavelength corresponds to the small parameter $\ep>0$, and the amplitude 
$a^\ep$ admits an asymptotic (in general not convergent) expansion, 
$$a^\ep \sim \sum_{n\geq0} \ep^n a_n(t,x), \quad \mbox{as } \ep\rightarrow0.$$
In particular, Lax uses such solutions to study the well-posedness of the 
Cauchy problem. 

The geometric information is contained in the \emph{eikonal equation} 
(of Hamilton-Jacobi type) determining the phase $\ph$ (plane waves
correspond to phases which are linear in $(t,x)$; spherical waves, to 
phases which are functions of $t$ and $|x|$ only). The amplitudes $a_n$ 
are solutions to \emph{transport equations} along the rays associated with 
the eikonal equation. The asymptotic expansion is valid for times
before rays focus.

\subsection{Weakly nonlinear geometric optics}

Trying to generalize Lax's approach to  
nonlinear systems, one immediately faces several problems. First, a family 
$\left(a(x)e^{i\ph(x)/\ep}\right)_{\ep>0}$ is not bounded in any
Sobolev space $H^s$ with  
$s>0$, so that the time of existence of the considered exact solution $u^\ep$ 
may shrink to zero as $\ep$ goes to zero. In the same spirit, the size of 
$u^\ep$ must be adapted, large enough to allow nonlinear features at
first order, but sufficiently small to prevent blow-up. Finally, the
main interest and difficulty of nonlinear models comes from
interactions: one hopes to incorporate in the asymptotic description
the self-interaction of a wave (leading for example to generation of
harmonics) as well as the possibility of (resonant) interaction of
several waves.

Thus, the extension to nonlinear systems goes through \emph{profiles}
$U_n(t,x,\theta)$ which are periodic or almost periodic w.r.t. 
$\theta\in\R^q$,
\begin{equation} \label{WNLexpansion}
u^\ep \sim \underline{u}(t,x) + \ep^m \sum_{n\geq0} \ep^n U_n
\left( t,x,\vp(t,x)/\ep \right),
\end{equation}
where $\underline{u}$ is a given groundstate, and
$\vp=(\ph_1,\dots,\ph_q)$ is a collection of phases.

The usual strategy of nonlinear geometric optics consists in:
1-defining a \emph{formal} solution, \ie solvable equations for
the profiles; truncating the series in~\eqref{WNLexpansion}
defines an approximate solution (a function approximately solution to
the equation); 2-showing that for any initial data close to the initial 
value of the approximate solution, an exact solution exists on a time
interval independent of $\ep$; 3(\emph{stability})-showing that
the exact solution is well approximated by the formal one.

In order to observe some nonlinear behavior, the magnitude $\ep^m$ of
oscillations is chosen so that cumulated effects of nonlinearities
become of the same order as the wave on the typical time $T$ of
propagation. A fixed $T$ (w.r.t. $\ep$) leads to \emph{nonlinear
  geometric optics}. Once rescaled, the system of partial differential
equations takes the form
\begin{equation} \label{syst}
L(t,x,u^\ep,\ep\partial) u^\ep = F(t,x,u^\ep),
\end{equation}
where the operator $L$ is in general a first order symmetric
hyperbolic system on some domain in $\R^{1+d}_{t,x}$,
\begin{equation} \label{op}
\begin{split}
L(t,x,u,\ep\partial)~u 
& = \ep\dt u + \sum_{j=1}^d A_j(t,x,u) 
\ep\partial_{x_j} u^\ep + L_0(t,x)~u \\
& =: L_1(t,x,u,\ep\partial)~u + L_0(t,x)~u.
\end{split}
\end{equation}
The smooth functions $F$, $B$ and $A_j$ take values in $\C^N$, in the
space of $N\times N$ matrices, and in the space of $N\times N$
symmetric matrices, respectively. Furthermore, we will distinguish 
the hyperbolic case, when $L_0=0$, and the dispersive case, when $L_0$ 
is skew-symmetric (the dissipative case, when
$L_0+L_0^\star \geq 0$, will be considered only in
Section~\ref{seclong}). This coefficient $L_0$ reflects the interaction 
between the wave and the material medium, so that, in the dispersive 
case, the group velocity (see~\eqref{vg}) depends on the frequency. The 
abovementioned exponent $m$ is then defined as follows. Let $J\geq2$ 
be the order of nonlinearities,
$$|\alpha| \leq J-2 \quad \Rightarrow \quad \partial^\alpha_{u
  \overline{u}} {A_j}_{|_{u=0}} = 0 , \qquad |\beta| \leq J-1 
\quad \Rightarrow \quad \partial^\beta_{u \overline{u}} {F}_{|_{u=0}} = 0.$$
The standard amplitude of weakly nonlinear geometric optics
corresponds to
\begin{equation} \label{WNLampl}
m = \frac{1}{J-1}
\end{equation}
(so that $m=1$ for quadratic semi- or quasilinear systems).

As an example, Euler equations for compressible and entropic gas
dynamics form a quasilinear hyperbolic system ($L_0=0$), 
and Maxwell-Bloch equations, a semilinear dispersive
one (describing the propagation of an electromagnetic wave $(E,B)$ in
a medium with polarization vector $P$; for a two-level quantum medium,
$N$ is the difference of the populations of the energy levels), which
reads
\begin{equation} \label{MB2}
\left\{
\begin{split}
& \dt B + \curl E = 0 , \\
& \dt E - \curl B = -\dt P , \\
& \dt N = E\cdot\dt P , \\
& \ep^2\dt^2 P + P = (N_0-N) E.
\end{split}
\right.
\end{equation}
It may be written in the form \eqref{syst} for the unknowns
$B,E,N,P,\ep\dt P$. 

\subsection{Profile equations} \label{secprofile}

The formal derivation of profile
equations is similar in any dimension, and for any number of phases. 
Plug expansion \eqref{WNLexpansion} into equation \eqref{syst}, 
and let the expansion of $Lu-F(u)$ vanish. This gives an infinite set 
of equations,

\begin{equation*} 
L(\dd\vp\cdot\dth)~U_0 = 0 ,
\end{equation*}
and for $n\geq0$,
\begin{equation} \label{BKWsyst}
L(\dd\vp\cdot\dth)~U_{n+1} + L_1(\partial)~U_n 
+ F_n(U_0,\dots,U_n) = 0  ,
\end{equation}
with $\dd\vp\cdot\dth = \sum_{k=1}^q \dd\ph_k \partial_{\theta_k}$.

The analysis is based on formal Fourier series in $\theta$, so that
for $U = \sum_\alpha U^\alpha e^{i\alpha\cdot\theta}$,
$$0 =  L(\dd\vp\cdot\dth)~U = \sum_\alpha
L(i\dd(\alpha\cdot\vp))~U^\alpha e^{i\alpha\cdot\theta} \Longleftrightarrow
U = \Pi~U = \sum_\alpha \pi_\alpha U^\alpha e^{i\alpha\cdot\theta} ,$$ 
where $\pi_\alpha$ is the orthogonal projector on the kernel of
$L(i\dd(\alpha\cdot\vp))$. In order to get non-trivial solutions, this
projector is assumed not to vanish at least for one multi-index
$\alpha$. This means that the function $\alpha\cdot\vp$ satisfies the
equation 
\begin{equation} \label{eikdet}
\det L(i\dd(\alpha\cdot\vp)) = 0,
\end{equation}
so that $\alpha\cdot\vp$ must satisfy one of the eikonal equations,
$$\dt(\alpha\cdot\vp) + \lambda_k(t,x,\dx(\alpha\cdot\vp)) = 0 .$$
Here, we denote by $\lambda_k(t,x,\xi)$, $1 \leq k \leq N$, the $N$ 
eigenvalues of the symmetric matrix 
$$\frac{1}{i} L(t,x,0,(0,i\xi)) = \sum_{j=1}^d
\xi_j A_j(t,x,0) + \frac{1}{i} L_0(t,x). $$
In particular, a single
initial phase $\ph^0(x)$ may generate $N$ eikonal phases $\ph_k(t,x)$ 
(except when initial profiles are prepared, \ie polarized on a single mode).

When the linear combination $\alpha\cdot\vp$ is not
trivial, \ie contains linearly independent $\ph_k$'s, equation~\eqref{eikdet}
expresses a \emph{resonance} between these phases. Such resonances
must involve at least three phases, and thus occur for systems of size
at least 3. 

Profile equations are then obtained recursively, 
splitting~\eqref{BKWsyst} thanks to the projectors $\Pi$ and $(1-\Pi)$,
\begin{subequations} \label{eqprof}
\begin{align}
& (1-\Pi) ~ U_0 = 0 , \qquad \mbox{ and for } n\geq0 , \label{polar} \\
& \Pi~L_1(\partial)~\Pi~U_n + \Pi~F_n = 0 , \label{transp} \\
& (1-\Pi) ~ U_{n+1} = - L(\dd\vp\cdot\dth)^{-1} [ L_1(\partial)~U_n +
F_n  ] . \label{ellipt}
\end{align}
\end{subequations}

Equation~\eqref{polar} is a \emph{polarization} constraint on the 
first profile $U_0$, and \eqref{transp} gives the evolution for the 
polarized part $\Pi~U_n$ of $U_n$ in the domain $\Omega \subset
\R^{1+d}$ considered. The first equation~\eqref{transp} (for the first
profile $U_0$) is nonlinear, while the others (for $n\geq1$) are
linear. The operator $L(\dd\vp\cdot\dth)^{-1}$ is  
formally defined by the symbol $L(i\dd(\alpha\cdot\vp))^{-1}$, the 
partial inverse of $L(i\dd(\alpha\cdot\vp))$ on its range. 

The description of the operator $\Pi~L_1~\Pi$ depends on the geometry 
of the characteristic variety,

$$\mathcal{C} = \{ (t,x,\tau,\xi) \in \Omega\times\R^{1+d} \mid \det
L(t,x,\tau,\xi) = 0 \},$$ 
viewed as a differentiable manifold. In the hyperbolic case ($L_0=0$), 
this characteristic variety is homogeneous in $(\tau,\xi)$: for the 
wave equation, it is the light cone $\mathcal{C}=\Omega\times
\{\tau^2=|\xi|^2\}$. In the dispersive case ($L_0=L_0^\star\neq0$), 
only few harmonics of an eikonal $(\tau,\xi)$ are eikonal in general: 
for Klein-Gordon equation, $\mathcal{C}=\Omega\times\{\tau^2=|\xi|^2+1\}$.

When for all $(t,x)\in\Omega$, the vector $\dd(\alpha\cdot\vp)(t,x)$
is a smooth point of a sheet $\tau=\lambda(\tau,\xi)$ of
$\mathcal{C}$, the principal part of $\pi_\alpha L_1(\partial)\pi_\alpha$ 
is simply the (scalar) \emph{transport operator} at the 
\emph{group velocity}, 
\begin{equation} \label{vg}
\vec{v}_\alpha(t,x) =
-\frac{\partial\lambda}{\partial_\xi}(t,x,\partial_x(\alpha\cdot\vp)(t,x));
\end{equation}
(for the wave equation, $\partial_\xi(\pm|\xi|) = \pm \xi/|\xi|$ does
not depend on $|\xi|$, whereas for Klein-Gordon equation,
$\partial_\xi(\pm\sqrt{|\xi|^2+1}) = \pm \xi~/\sqrt{|\xi|^2+1}$ does).  
Precisely (\cite{Lax57}, \cite{DR97a}),
$$\pi_\alpha L_1(\partial)\pi_\alpha = \dt + \vec{v}_\alpha\cdot\dx +
\pi_\alpha(L_1(\partial)\pi_\alpha).$$ 

This reduction clearly explains the ``light rays'' picture; it is
however not necessary in the justification of the
asymptotics~\eqref{WNLexpansion}; as an example, Joly, M\'etivier and
Rauch perform this justification \cite{JMR94} in the case of an
eigenvalue $\lambda$ which changes multiplicity (crossing in the
characteristic variety), which gives a model for conical refraction
of light. Moreover, this description is clarified by Lannes 
\cite{Lan98}, in the case of isolated singular points of a dispersive
$\mathcal{C}$: the characteristic variety of $\pi_\alpha
L_1(\partial)\pi_\alpha$ is then the tangent cone to $\mathcal{C}$. A
systematic approach to such ``algebraic lemmas'' is given by Texier 
\cite{Tex04} \via perturbation theory of complex eigenvalues under a
slightly different smoothness assumption (found by M\'etivier and
Rauch \cite{MR03} to be actually the same in the hyperbolic case). 
This approach also generalizes the one of Lannes, and applies to 
``all order algebraic lemmas'' (see Section~\ref{secmorescales}).

\subsection{Rigorous results} 

One easily suspects, in view of the term
$L(\dd\vp\cdot\dth)^{-1}$ in \eqref{ellipt}, that resonances may 
cause troubles. That is the reason why monophase expansions (with
profiles polarized on a single mode) have been studied first. 

In \cite{CB69}, Choquet-Bruhat constructs formal such expansions 
for quasilinear systems, and applies these to computations
on the model of relativistic perfect fluids. Single-phase weakly
nonlinear geometric optics is justified in space dimension $d$, for 
semilinear systems, by Joly and Rauch \cite{JR92} in the
non-dispersive case, by Donnat and Rauch \cite{DR97a}, \cite{DR97b} in
the dispersive case, and by Gu\`es \cite{Gue93} for quasilinear
non-dispersive systems. Their strategy consists in constructing an 
approximate solution $u^\ep_{\rm app}$, sum of a high number of smooth
profiles (or with infinitely accurate asymptotics, \via Borel's
summation technique), so as to absorb the oscillations of the residual $L
u^\ep_{\rm app} - F(u^\ep_{\rm app})$. Existence of exact solutions
$u^\ep$ and stability are then obtained \via energy estimates on
$u^\ep-u^\ep_{\rm app}$. 

The formal study of propagation and interactions of oscillations for
quasilinear systems, particularly for fluid dynamics, is done in the
works of Hunter and Keller \cite{HK83} in the non interacting case, of
Majda and Rosales \cite{MR84} in the one dimensional resonant case, of
Hunter, Majda and Rosales \cite{HMR86} and Pego \cite{Peg88} in the
multidimensional resonant case. In these works, a finiteness
assumption is made on the number of (directions of) eikonal
phases. But in \cite{JR91}, Joly and Rauch show that resonances may
generate dense oscillations in the characteristic variety. 

This finiteness hypothesis is shown unnecessary in \cite{JMR93a},
where Joly, M\'etivier and Rauch justify rigorously the weakly
nonlinear geometric optics approximation in space dimension 1, for
both semi- and quasilinear systems (before shock formation), under a
weaker transversality assumption between real combinations $s$ of phases
and the propagation fields $X_k$ ($X_k s \equiv 0$, or $X_k s \neq
0$ almost everywhere on the domain $\Omega$). They use the method
initiated by Joly 
\cite{Jol83} for semilinear multidimensional systems with constant
coefficients and linear phases, estimating (by stationary and
non-stationary phase method) the difference $u_\ep^k - \ep^m 
U^k_0(t,x,\vp/\ep)$ between the $k$-th Picard iterates obtained from the
resolution of~\eqref{syst} and of the first profile equation,
respectively. This leads to a first order approximation,
$$u^\ep(t,x) = \ep^m U_0(t,x,\vp(t,x)/\ep) + o(\ep^m),$$
which is valid in $L^p$ for all $p<\infty$. It holds in $L^\infty$
when transversality is imposed (almost everywhere) along integral
curves of each propagation field $X_k$. Gu\`es shows \cite{Gue95a} that
on the contrary, when 
weak resonances are present (a linear combination of phases is eikonal
on a set with one dimensional positive measure, but is a.e. not eikonal
on $\Omega$), they may lead to the creation of stationary waves
preventing from $L^\infty$ approximation. 

Concerning results in space dimension one, we have the earlier works
(for linear phases) of Kalyakin \cite{Kal89}, and also of Tartar
\emph{et al.} \cite{Tar81}, \cite{Tar83}, \cite{McLPT85}, by the compensated
compactness method applied to Young measures --links between nonlinear
geometric optics and compensated compactness are very clearly
enlighted in \cite{MS98} and \cite{JMR95b}.

Furthermore, still in space dimension one, several papers deal with
(global) weak solutions (with bounded variations) of conservation
laws. The formal study is due to Majda and Rosales \cite{MR84};
DiPerna and Majda obtain $L^\infty_t L^1_x$ global in time asymptotics
--even when shocks occur--, using $L^1$ stability and BV decay
properties of solutions (constructed by Glimm's scheme) to
conservation laws, for BV initial data with compact support. 
This result is extended by Schochet \cite{Sch94a} to
periodic BV data using WKB expansions, and then by Cheverry
\cite{Che96b}, \cite{Che97} to general BV data. Sabl\'e-Tougeron
\cite{Sab96} treats the initial-boundary value problem, for linear
non-interacting phases, and data with compact support.

Now, the multidimensional case for semi- or quasilinear interacting
waves is treated by Joly, M\'etivier and Rauch \cite{JMR95a}. In space
dimension greater or equal to 2, the new phenomenon is the
\emph{focusing} of oscillations, corresponding to singularities of
phases (see Section~\ref{seccaustics}). It may lead to blow-up, and to
the ill-posedness of the Cauchy problem. Joly, M\'etivier and Rauch
distinguish ``direct focusing'' of principal phases and ``hidden
focusing'' resulting from several interactions. They define a
\emph{coherence} criterium so as to separate clearly eikonal and
non-eikonal phases. 

\begin{defn}
A real vector space $\Phi \subset \mathcal{C}^\infty(\Omega)$ is
$L$-coherent when for all $\ph\in\Phi\nz$, one of the following
conditions holds: \\ 
{\rm (i)} $\det L(t,x,\dd\ph(t,x)) = 0$ and $\dd\ph(t,x)\neq0$ for
every $(t,x)\in\Omega$, \\
{\rm (ii)} $\det L(t,x,\dd\ph(t,x))$ does not vanish on $\Omega$.
\end{defn}

The typical example of $L$-coherence is for linear phases and constant
coefficient operator $L$. Coherence is not easy to check on a given
set of phases, but it
allows to prove (locally, on truncated cones $\Omega_T= \{(t,x) \mid 0
\leq t \leq T~,~|x|+t/\delta\leq\rho\}$, for small $\delta,\rho>0$),
for a symmetric hyperbolic system with constant multiplicity, the
following result.

\begin{theo} Consider initial phases
  $\vp^{~0}=(\ph^0_1,\dots,\ph^0_q)$ such that
  $\dd \ph^0_k$ does not vanish on $\Omega_0$. Consider eikonal phases 
  $\vp=(\ph_1,\dots,\ph_q)$ generated by these initial
  phases ($\forall k, ~ \ph_k(0,x)=\ph^0_k(x)$), and assume that their
  linear span $\Phi$ is $L$-coherent, 
  and contains a timelike phase $\ph_0$: $\ph_0(0,x)=0$, and
  $\dt\ph_0$ does not vanish. Then, for any initial data
\begin{equation} \label{inidata}
h^\ep(x) = \ep^m\mathcal{H}^\ep(x,\vp^{~0}/\ep),
\end{equation}
with $(\mathcal{H}^\ep)_{0<\ep\leq1}$ a bounded family in the Sobolev
space $H^s(\overline{\Omega_0}\times\T^q)$ for some $s>1+(d+q)/2$, we
have \\
{\rm (i)} There is $T>0$ such that for all $\ep \in ]0,1]$, the Cauchy
problem~\eqref{syst}, \eqref{inidata} has a unique solution
$u^\ep\in\mathcal{C}^1(\Omega_T)$. \\
{\rm (ii)} There is a profile $\mathcal{U}(t,x,\tau,\theta) \in
\mathcal{C}^1(\Omega_T\times\R\times\T^q)$, almost periodic in
$(\tau,\theta)$ (and determined by the profile
equations~\eqref{polar}, \eqref{transp}) such that
$$u^\ep(t,x) - \ep^m \mathcal{U}(t,x,\ph_0/\ep,\vp/\ep)
= o(\ep^m) \mbox{ in } L^\infty.$$
\end{theo}

The new idea in their proof consists in looking for exact solutions of
the form
$$u^\ep(t,x) = \ep^m \mathcal{U}^\ep(t,x,\vp/\ep).$$
It is then sufficient (for $u^\ep$ to be a solution to~\eqref{syst})
that the profile $\mathcal{U}^\ep$ be solution to the \emph{singular
  system}
$$L_1(\partial_{t,x})~\mathcal{U}^\ep + \frac{1}{\ep} \left(
  \sum_{k=1}^q L_1(\dd\ph_k) \partial_{\theta_k} + L_0
\right)~\mathcal{U}^\ep = F^\ep(t,x,\mathcal{U}^\ep),$$
and $L$-coherence precisely allows energy estimates in $H^s$.

The result of Joly, M\'etivier and Rauch
on singular systems generalizes the ones of Klainerman and Majda
\cite{KM81}, \cite{KM82} and Schochet \cite{Sch94b}. It is worth
noting that the analysis of such singular systems appears 
in many multi-scale problems. In particular, for rapidly rotating fluids
in oceanology, meteorology, and Magneto-Hydro-Dynamics; see
\cite{BMN96}, \cite{Gre97a}, \cite{Gal98}, \cite{MS01}; several 
applications are described by Klein \cite{Kle00}. This method also
provides an efficient tool for numerics: see the works of Colin and
Nkonga about propagation of wavetrains \cite{CNk04} and pulses (see
section \ref{secothers}) \cite{CN05} in optical media, and Colin and
Torri \cite{CTp} about pulse propagation over diffractive scales (see
section \ref{seclongtime}).  

A work close to the method exposed here may be found in the paper
\cite{CC} by Colin and Colin, studying rigorously and numerically
Raman scattering for a semi-classical quasilinear Zhakarov system from
plasma physics in the weakly nonlinear geometric optics regime, with
three-wave resonances. There, the system lacks hyperbolicity, which
is compensated by dispersion.  

It is in fact possible to achieve infinitely accurate asymptotics for
these multidimensional nonlinear interacting waves with nonlinear
phases, by adding to the coherence assumption on phases a generically
satisfied ``not too small divisors'' hypothesis (or Diophantine
hypothesis on wave vectors, in the case of linear phases), see
\cite{JMR93b}.

\section{Other kinds of profiles} \label{secothers}

As we have seen, the natural dependence of the profiles
$U(t,x,\theta)$ on the fast variable $\theta\in\R^q$ is periodic,
quasi- or almost periodic. This allows to define an oscillating
spectrum \cite{Jol83}, \cite{JMR94} which is localized on the
frequencies given by the Fourier transform of $U(t,x,\cdot)$. Now,
WKB aymptotics with profiles having other behaviors w.r.t. $\theta$
may be relevant, depending on the context. The formal
computations are in general very similar to the ones of usual
geometric optics, but at least, the functional tools (such as average
operators) must be re-defined, and interactions may take a different form.

In order to model the optical Raman scattering, for which light is
emitted in a continuum of directions,
Lannes \cite{Lan98} introduces profiles with continuous spectra,
and gives the analogue of usual rigorous weakly nonlinear geometric
optics in this context, with a precise analysis of resonances. 

In fact, this formalism also includes (see Barrailh and Lannes
\cite{BL02}) the one of ``ultrashort pulses''
(from laser physics) considered by Alterman and Rauch \cite{AR02},
where profiles have a compact support in $\theta$ (see also 
sections~\ref{seccaustics} and~\ref{secmorescales} about focusing and
diffraction of pulses). 

Sometimes, profiles may also have
different limits as the variable $\theta\in\R$ goes to $\pm\infty$, in
order to match boundary conditions (see Section~\ref{bdylayers}), or
to describe transitions like the ones of \cite{Dum04}, between light and
shadow (see Section~\ref{secmorescales}); this is also the case of the
solitary waves in \cite{Gue95a}. 

\section{Caustics} \label{seccaustics}

\begin{figure}[htbp]
$$
\centerline{\hbox{\input{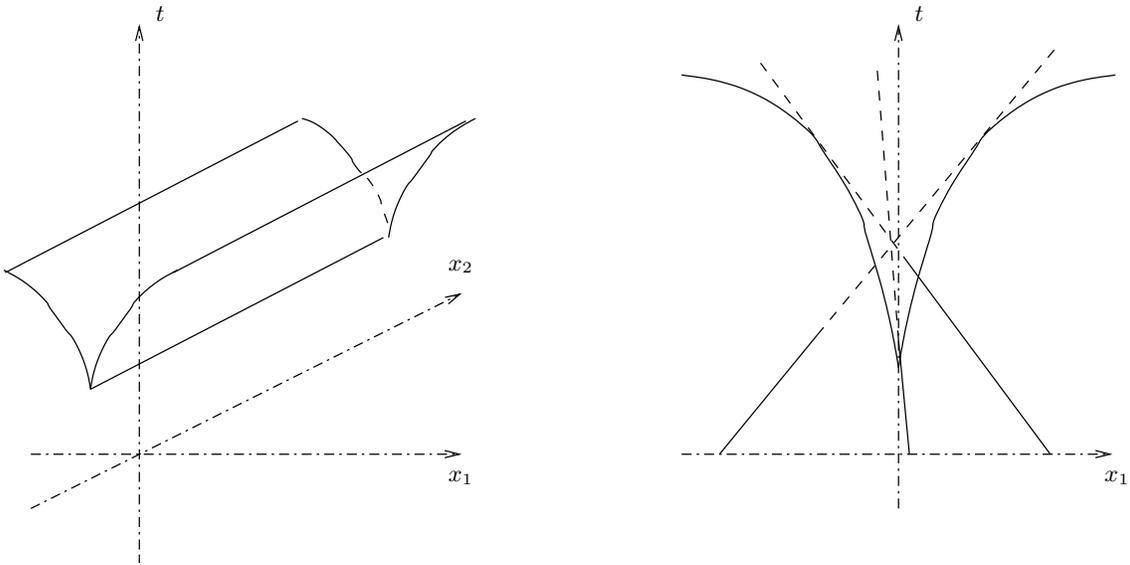}}}
$$
\caption{The cusp. Through each $(t,x)$ point pass only one ray
  before the caustic, and three rays beyond the caustic.} 
\label{cusp}
\end{figure}

In space dimension greater or equal to 2, singularities appear, even
in the case of linear geometric optics, in the resolution of the
eikonal equation: when rays have an envelope (called a caustic; the
example of the cusp $\{(t,x_1) \in [0,\infty[\times\R \mid t^{2/3} =
x_1^{2/3}+1\}$ from the wave equation in space dimension 2 and
$\ph(0,x)=x_2+x_1^2$ is shown on Figure~\ref{cusp}), the second derivative
of the phase $\ph$ becomes singular. In this case, the amplitude is
also singular \cite{Lud66}.

In order to solve the eikonal equation globally in time, one considers 
the Lagrangian manifold folliated by the Hamiltonian flow associated
with $L(t,x,\tau,\xi)$, starting from points
$(0,x,0,\dd\ph(0,x))$. Caustic crossing induces a phase shift
determined by Maslov's index, and generates new phases (see Duistermaat
\cite{Dui74}; one phase before the cusp of Figure~\ref{cusp} corresponds to
three phases beyond the cusp).

Following these ideas, Joly, M\'etivier and Rauch have studied the
caustic crossing for semilinear geometric optics. In \cite{JMR95c},
they show that, for superlinear nonlinearities and focusing at a
point, oscillations may lead to explosion, whereas in the case of
dissipative equations (for which exact solutions are globally
defined), they may be absorbed (\ie only a non-oscillating term
remains after reflection on the caustic). In \cite{JMR00a}, 
for dissipative equations, they extend
this result to general caustics. They exhibit a critical exponent
$p_c$, defined by the geometry of the caustic so that, if the
nonlinearity is stronger than $|u|^{p_c}$ at infinity, absorption
occurs, and else, oscillations persist. In \cite{JMR96}, for uniformly
Lipschitzean nonlinearities (for which, again, exact solutions are globally
defined), they show that oscillations go through the caustic. 

More precisely, to a single-phase expansion before the caustic
(case considered here just for notational simplicity; several initial
phases are allowed),
$$u^\ep(t,x) = \ep^m \mathcal{U}(t,x,\ph/\ep) 
+ o(1) , $$
corresponds a multiphase expansion beyond the caustic,
$$u^\ep(t,x) = \ep^m \underline{u}(t,x) + \ep^m \sum_j
\mathcal{U}(t,x,\ph_j/\ep) + o(1) , $$ 
with phases $\ph_j$ defined as in the linear case, matching $\ph$ on
the caustic. In particular, this provides examples of asymptotics with
non-coherent phases. This approximation is in $L^2$ in \cite{JMR95b} and
\cite{JMR96}, and in all $L^p$'s, $p<\infty$, in \cite{JMR00a}. It is
based on the oscillating integral representation 
\begin{equation*}
\begin{split}
u^\ep(t,x) 
& = \ep^m \underline{u}(t,x) + \ep^m \int e^{i\phi(t,x,y,\xi)/\ep}
\mathcal{A}(t,x,\ph/\ep)~\dd y~\dd\xi + o(1) \\ 
& =: \ep^m \underline{u}(t,x) + \ep^m I^\ep(\mathcal{A}) + o(1) . 
\end{split}
\end{equation*}
Outside the caustic, expanding a smooth
$\mathcal{A}(t,x,\theta)$ in Fourier series in $\theta$ and performing 
the usual stationary phase asymptotics (the assumption that the
phases $\ph_j$ are non-resonant is needed) gives 
$$I^\ep(\mathcal{A}) \sim \sum_j J_j^\ep(\mathcal{A}) ,$$ 
where the sum corresponds to the several pre-images of critical
points of the phase, and $J_j^\ep(\mathcal{A})(t,x)$ is the evaluation at
$(t,x,\ph_j/\ep)$ of some profile $\mathcal{U}_j(t,x,\theta)$ obtained
from $\mathcal{A}$ \via a Hilbert transform.

Here, the amplitude $\mathcal{A}$ belongs to $L^p$, so that the
profiles $\mathcal{U}_j$ are in fact better defined as ``weak
profiles'' (in particular, $\mathcal{U}_j(t,x,\ph_j/\ep)$ has an
asymptotic sense only). They are weak limits of $u^\ep$, using
$J_j^\ep(\mathcal{B})$ as test function, for smooth
$\mathcal{B}(t,x,\theta)$,
$$\forall \mathcal{B}(t,x,\theta), \quad \int \overline{u^\ep}
J_j^\ep(\mathcal{B}) \dd x \tendlorsque{\ep}{0} \int
\overline{\mathcal{U}_j} \mathcal{B} \dd x \dd\theta.$$ 
Profile equations are then obtained as weak limits of the original
equation, using a nice ``no propagated oscillations'' lemma
(\cite{JMR00a}, section 5.2) for nonlinearities. Strong convergence
then follows from energy estimates and $L^p$ estimates on oscillatory
integrals. 

These results are refined (in $L^\infty$ instead of $L^p$) and
extended to conservative systems by Carles, for simpler
geometries. Essentially, Carles considers focusing on a point for the
wave equation, 
\begin{equation} \label{W}
\square u^\ep + a |\dt u^\ep|^{p-1}\dt u^\ep = 0, \quad p>1, ~ a\in\C,
\end{equation}
($a>0$, $a<0$, $a\in i\R$ corresponds to the dissipative,
accretive and conservative case, respectively) with radially symmetric
initial data, or the semiclassical nonlinear Schr\"odinger equation (NLS),
\begin{equation} \label{NLS}
i\ep\dt u^\ep + \frac{1}{2} \ep^2 \Delta u^\ep = \ep^\alpha
|u^\ep|^{2\sigma}u^\ep, \quad \alpha\geq1, ~ \sigma>0, 
\end{equation}
with initial data of the form
$$u^\ep_{|_{t=0}} = f(x) e^{-i|x|^2/2\ep}.$$

Hunter and Keller \cite{HK87} give a formal classification of the
qualitative properties of weakly nonlinear geometric optics for
\eqref{W}, separating linear and nonlinear propagation, and linear and
nonlinear effects of the caustic. Carles has rigorously explained this
classification. In \cite{Car98}, for $1<p<2$, he improves the results
of \cite{JMR00a} (\cite{JMR95c} corresponds to $a>0$, $p>2$) with an 
$L^\infty$ description of radial waves in $\R^3$. He shows that the
profiles are really singular, and that new amplitudes (of size
$\ep^{1-p}$) must be added to the one in \cite{HK87} near the caustic. 

In a series of papers \cite{CR02}, \cite{CR04a}, \cite{CR04b}, Carles
and Rauch study the focusing of spherical pulse solutions to \eqref{W}
in space dimension 3, for various powers $p$ and sizes ($\ep^m$) of 
initial data. In particular, they show how pulses get wider
after reflection on the caustic; they also treat the analogue of
\cite{JMR95c} ($a\in\R$) and find the same results (dissipation or
growth) as for wavetrains. Carles and Lannes generalize these results
\cite{CLp} to dispersive semilinear systems such as Klein-Gordon and
Hartree's equations. Again for \eqref{W}, they obtain \cite{CL03} 
the focusing (at $t=1$) of a pulse with ``random phase shift'' 
$\ln ((1-t)/\ep)$ --in the conservative case, for nonlinear
propagation and a nonlinear caustic. 

Concerning NLS \eqref{NLS}, Carles investigates all the behaviors with
different powers $\alpha\geq1$ and $\sigma>0$. In \cite{Car01}, he
shows the possibility of ``random phase shift''. In \cite{Car00a}
($\alpha=d\sigma>1$, in space dimension $d$), he notices that
nonlinear effects take place near the focal point only; thus, in
\cite{CFKG03}, \cite{CKp}, he answers with his co-authors the converse
question ``for which kind of initial data can we get nonlinear effects
at first order?''. He also considers \cite{Car00b} focusing on a line
instead of a point. 

Other extensions of \cite{Car00a} concern the description of
Bose-Einstein condensates, \via the addition to \eqref{NLS} of a
polynomial potential \cite{Car03}, \cite{CM04}, \cite{CN04}, or 
Hartree-type nonlinearities $((1/|x|) \star |u|^2) u$ \cite{CMSp} 
for the Schr\"odinger-Poisson model. 

The method of Carles consists in a precise control of the solution
outside the caustic \via WKB expansions and Gagliardo-Nirenberg
estimates, then rescaling the problem around the focal point so as to
use nonlinear scattering arguments (with short-range or long-range
interpretation, according to the power of the nonlinearity).

We may mention that caustics are subject to many numerical works, such
as the ones of Benamou {\it et al.}, \cite{BS00}, \cite{BH02},
\cite{Ben03}, \cite{BLSS03}, \cite{BLSS04}, based on a Eulerian --or
ray tracing-- approach. See also the review of Gosse \cite{Gosp}
--about (NLS), mostly with kinetic methods-- and references therein.

\section{Boundary problems}

\subsection{Reflection on a boundary}

The formal study of reflections of weakly nonlinear waves is done by
Artola and Majda in \cite{MA88}. In the case of transverse reflection,
rigorous results are due to Chikhi \cite{Chi91} for equations with two
speeds of propagation, and to Williams \cite{Wil93} for general
resonant interactions. 

The first order tangential reflection corresponds to diffractive
points in the cotangent bundle of the boundary, and leads to the
formation of a shadow zone (in fact, a caustic forms in the obstacle,
so that phases become singular at the tangency point). The linear case
is treated by Friedlander \cite{Fri76} and Melrose \cite{Mel75a},
\cite{Mel75b}. Weakly nonlinear geometric optics is justified, at
first order, with nonlinear incident and reflected phases, by Cheverry
\cite{Che96a} for globally Lipschitz nonlinearities (constructing smooth
profiles, solutions to singular ODE's along the broken flow of the
rays), and by Dumas \cite{Dum02} for a dissipative nonlinear
Klein-Gordon equation (using weak profiles as described in
section~\ref{seccaustics}). 

Higher order tangency is studied by Williams \cite{Wil93},
\cite{Wil00}, in the case of a simpler geometry (flat boundary,
constant coefficient operators and linear phases). He obtains
$L^\infty$ asymptotics at all orders, incorporating in the description
boundary layer profiles (see section~\ref{bdylayers}) involving a
third scale $\rep$: profile equations are then of Schr\"odinger type 
(see section~\ref{secmorescales}).

\subsection{Boundary layers} \label{bdylayers}

Boundary layers are a usual feature in the study on partial
differential equations on domains with boundaries: they stem from
large variations of the solution in a small region near the
boundary. They often appear in vanishing viscosity models, since the
boundary conditions are not the same for the viscous (parabolic)
and inviscid (hyperbolic) equations.

Gu\`es \cite{Gues95b} describes such boundary layers for viscous
perturbations ($-\ep^2\mathcal{E}$, with $\mathcal{E}(t,x,\partial)$ a
second order elliptic operator) of semilinear mixed problems
(associated with a linear hyperbolic operator $L(t,x,\partial)$). He
constructs exact solutions with infinitely accurate WKB asymptotics of
the form
$$u^\ep(t,x) \sim \sum_{n\geq0} \sqrt{\ep}^{~n}~\mathcal{U}_n \left(
  t,x,\ph/\rep,\ph/\ep \right) , $$
$$\mathcal{U}_n \left( t,x,\theta,z \right) =
\underline{\mathcal{U}}_n(t,x) + \mathcal{V}_n(t,x,\theta) +
\mathcal{W}_n(t,x,z) + \chi_n(t,x,\theta,z) ,$$
where the phase $\ph$ is transverse to the boundary, and functions of
$\theta, z \geq 0$ decay rapidly at $+\infty$ (so that
$\underline{\mathcal{U}}_n = \Lim{\theta,z}{+\infty}\mathcal{U}_n$).
Profile equations are analogue to the ones of geometric optics (but of
parabolic type), with polarization conditions. When the boundary is
non-characteristic for the hyperbolic operator $L$, the dependence
on $\theta=\ph/\rep$ is not necesseray, and the boundary layer has size
$\ep$; in the characteristic case, the boundary layer is wider, of
size $\rep$. Grenier \cite{Gre97b} considers the quasilinear case with
characteristic boundary, and Grenier and Gu\`es \cite{GG98}, the
non-characteristic quasilinear case. Extensions may be found
in the work of Sueur for semilinear systems \cite{Suea}, \cite{Sueb}, 
\cite{Sue05}, or quasilinear systems \cite{Suec}, \cite{Sued}. 

When oscillations are present in the interior domain, they may
interact with boundary layers. This occurs for example in the study of
glancing oscillations by Williams \cite{Wil00}. It also occurs for
rotating fluid models (see the books of Greenspan \cite{Gre80} and
Pedlovsky \cite{Ped79} about oceanography and meteorology, and Dormy 
\cite{Dor97} about magneto-hydro-dynamics), which in the
incompressible case typically have the form
$$\dt u^\ep + u^\ep\cdot\nabla u^\ep + \frac{1}{\ep} R u^\ep - \ep
\Delta u^\ep = -\nabla p^\ep ,$$
where $p^\ep$ is the pressure, $u^\ep$ is the divergence free velocity
field, and the operator $R$ is skew-symmetric (for example, $R u =
e \times u$ with a fixed vector $e$ for Coriolis effect). Initial data
are well-prepared when they correspond to admissible data for the
limit problem. Otherwise, fast time oscillations appear in the interior
domain, and are taken into account in the profiles through a
dependence on the variable $t/\ep$. For such studies in the
well-prepared case, see \cite{GM97}; in the ill-prepared case, see
\cite{Mas00}, \cite{CDGG02}, \cite{BDGV04}, \cite{Ger03a};
G\'erard-Varet emphasizes the link with geometric optics in
\cite{Ger03b} and \cite{Ger05}. 

This technique is also used for micromagnetism. Carbou, Fabrie and
Gu\`es \cite{CFG02} describe the magnetization of a ferromagnetic
medium, given by Landau-Lifshitz equation \cite{LL69}, as the exchange 
coefficient goes to zero; Sanchez \cite{San02} studies the diffraction 
of an electromagnetic wave by a conducting domain, as the conductivity 
goes to infinity.

\subsection{Shocks} \label{secshocks}

Shock waves entering the framework of weakly nonlinear geometric 
optics are weak shocks, since their amplitude goes to zero with
$\ep$. Such an oscillatory shock is a perturbation of a
non-oscillating shock $(\underline{u}^\pm,\underline{\psi})$, and is 
given (in space dimension $d$, with $x=(x',x_d)$) by a $\mathcal{C}^1$ 
hypersurface,
$$x_d = \psi^\ep(t,x') \sim \underline{\psi}(t,x) + \ep^m \sum_{n\geq0} 
\ep^n \psi_n(t,x',\vp^{~0}(t,x')/\ep),$$
and on each side $x_d \gtrless \psi^\ep(t,x')$, a $\mathcal{C}^1$ 
solution $u^\ep_\pm$ to a system of conservation laws, 
$$u^\ep_\pm(t,x) \sim \underline{u}^\pm(t,x) + \ep^m \sum_{n\geq0}
\ep^n \mathcal{U}^\pm_n(t,x,\vp^{~\pm}(t,x)/\ep),$$
together with the usual Rankine-Hugoniot jump conditions. Of course, 
the unknowns are $u^\ep_\pm$ and $\psi^\ep$, so that this is a 
free-boundary problem.

A formal description is due to Artola and Majda \cite{MA88}. Rigorous 
results in space dimension one are given by Corli \cite{Cor95} in a 
multiphase context, with curved phases, using the same ``Picard iterates" 
method as in \cite{JMR93a} (and the same transversality assumptions on 
phases). Williams \cite{Wil99} obtains, in the spirit of \cite{JMR93b}, 
multidimensional infinitely accurate asymptotics with smooth profiles, 
around a planar shock ($\underline{u}^\pm=cst$, $\underline{\psi}=\sigma 
t$), which is uniformly stable in the sense of Kreiss. This is done
for linear phases (with a Diophantine assumption on their gradients), 
whose restrictions to the shock surface are entire multiples of a
single phase $\ph^0$.
 
We also mention the work of Gu\`es, M\'etivier, Williams and Zumbrun on 
the existence and stability of multidimensional viscous shocks 
\cite{GMWZ04a}. These (non-oscillating) shocks are not weak ones (their 
size does not depend on a small parameter), and the system of conservation 
laws is perturbed by a viscous term $-\ep\Delta u^\ep$, so that the 
solution $u^\ep$ is smooth. The question is ``are the shock solutions of 
the inviscid system approximated by solutions of the viscous system?". 
The answer is ``yes" \cite{GMWZ04b}, under the assumption of the existence 
of a spectrally stable wave profile $U(x',(x_d-\psi)/\ep)$ solution to the 
viscous system such that $U(x',z) \tendlorsque{z}{\pm\infty} u^\pm$. 
These authors also study \cite{GMWZ05a}, \cite{GMWZ05b} the long time 
stability of planar shocks, using WKB asymptotics for the viscous problem. 

\section{Three (and more) scale expansions: diffractive optics}
\label{secmorescales} 

So far, we have essentially presented problems where only two scales appear 
(at least at first sight). But more scales may be present (remember how a 
third scale appears with boundary layers in section~\ref{bdylayers}). In 
this section, we shall see that a third scale leads to supplement the 
transport equations (with finite propagation speed) for the profiles
with a Schr\"odinger-type equation (with infinite propagation speed!)  
taking diffraction into account; this is the usual ``paraxial
approximation'' of laser physics in Kerr media. 

\subsection{Some three-scale problems}

In Donnat's thesis \cite{Don94}, a model of ``light cigars" is described 
for laser propagation. These are modulated oscillating wavetrains with 
frequency $1/\ep$, with anisotropic ``support", \ie of typical size 1
along a direction, and $\rep$ along another. Such profiles are used by 
Boucher\`es, Colin, Nkonga, Texier and Bourgeade \cite{BCNTB04} to
describe the interaction of a laser with a gas, giving an
approximation of Maxwell-Bloch equations \eqref{MB2} by a nonlinear
Schr\"odinger-Bloch system.  

``Singular rays" of weakly nonlinear geometric optics may also be
analysed with three-scale asymptotics. They are hypersurfaces
(defined as $\psi(t,x)=0$ for an unknown function $\psi$) constituted
by rays meeting an obstacle tangentially or at a corner. They model
transitions between light and shadow (or sound and silence, in the
case of acoustic waves). Hunter \cite{Hun88a} introduces formal WKB
expansions, 
\begin{equation} \label{singray}
u^\ep(t,x) \sim \ep^m \sum_{n\geq0}\rep^{~n} U_n(t,x,\psi/\rep,\vp/\ep),
\end{equation}
with profiles periodic w.r.t. $\theta=\ph/\ep$, and with limits as
$\eta=\psi/\rep \rightarrow \pm\infty$ (such a $u^\ep$ may be
understood as an smooth approximation of a contact discontinuity, \ie
a characteristic shock along $\psi=0$). With an approach close to the
one of \cite{DJMR96}, Dumas \cite{Dum04} validates this asymptotics
rigorously, introducing new coherence assumptions on the phases
$(\tilde{\psi},\vp)$, where $\tilde{\psi}=(t,\psi)$ --it turns out
that an intermediate phase $\psi_0=t$ is in general necessary for the
profile equations to be well-posed. Profile equations are then
transport equations in variables $(t,x)$, and Schr\"odinger equations
with time $\tau=t/\rep$ and space variables $x$.

Three-scale asymptotics of the form \eqref{singray} also describe
oscillations according to perturbed phases (in a more restricted
sense than in \cite{Gre98}). For initial data such as 
$$u^\ep_{|_{t=0}} = \ep^m \mathcal{U}^0 (x,\ph_\ep^0/\ep) \mbox{ with
} \ph_\ep^0 = \ph^0 + \rep \psi^0, $$
the solution to the Cauchy problem cannot in general be approximated 
in $L^\infty$ by a two-scale geometric optics description, but
profiles as above, satisfying a NLS equation on a torus, achieve such
an approximation --see Dumas \cite{Dum03a}. 

Now, we turn to the problem of long-time validity for geometric
optics, in which the third scale is not present initially, but
naturally comes up. This is the context where first rigorous 
three-scale asymptotics were proven in nonlinear optics.

\subsection{Long-time behavior} \label{seclongtime}

\emph{A priori}, asymptotics from section \ref{secfirst} are valid
only on some fixed time interval $[0,t_0]$. In fact, when the profiles
are globally defined, uniformly bounded in space-time, and have $H^s$ 
norms with polynomial growth in time, Lannes and Rauch prove
\cite{LR01} that weakly nonlinear geometric optics is valid up to
times $t\sim C\ln(1/\ep)$ --in the linear case, see the study of
propagation on this so-called Ehrenfest time for semi-classical
Schr\"odinger equations by De  Bi\`evre and Robert \cite{DBR03}.

\paragraph{Diffractive optics.} However, in order to describe
oscillatory waves with frequency $1/\ep$ 
on larger propagation scales, one must adapt the geometric optics
approach (even in the linear case). First, the cumulated effects of
nonlinearities over longer 
times leads to consider smaller amplitudes. The first natural
long-time scale is then the \emph{diffractive time scale} $1/\ep$, and
formula \eqref{WNLampl} is replaced with
\begin{equation} \label{difframpl}
m = \frac{2}{J-1}. 
\end{equation}
Next, new variables are introduced so as to capture evolution up to
$t\sim1/\ep$, $X=\ep x$ and the slow time $T=\ep t$.

Note that a solution to the eikonal equation is in general not smooth
globally in time, unless the operator $L$ has constant coefficients
and the phase is linear, of the form $\ph(t,x)=\beta\cdot x-\omega t$.

Formal three-scale asymptotics for solutions $u^\ep$ to \eqref{syst}
are given by Donnat \cite{Don94},
\begin{equation} \label{profdiff}
u^\ep(t,x) \sim \ep^m \sum_{n\geq0} \ep^n U_n
(T,X,t,x,\theta)_{|{(T,X,\theta)=(\ep t,\ep x, (\beta\cdot x-\omega
    t)/\ep)}}. 
\end{equation}
The natural dependence of the profiles $U_n$ (periodic in $\theta$) 
in the variable $X$ (and $x$) is of Sobolev type, and continuous (thus
bounded) in $T \in [0,T_0]$, whereas for consistency of the Ansatz,
\ie $|\ep U_{n+1}| \ll |U_n|$, \emph{sublinear growth} in $t=T/\ep$ is 
required,
$$\forall n\geq1, \quad \frac{1}{t} \|U_n\|_\infty
\tendlorsque{t}{\infty} 0.$$ 

Donnat, Joly, M\'etivier and Rauch \cite{DJMR96} give a rigorous
justification of this monophase approximation, for semilinear
hyperbolic (non-dispersive) systems,
and profiles with no non-oscillating part (assuming that the Taylor
expansion of the nonlinearities contains odd parts only). The
Schr\"odinger profile equation (see \eqref{eqprofdiffRNLS}) is interpreted
as a diffractive correction to geometric optics on long times, by
analogy with Fresnel's diffraction.

In order to get faster numerics for nonlinear optics, Colin, Gallice
and Laurioux (Barrailh) \cite{CGL05} introduce intermediate models between
geometric and diffractive optics, keeping a dependence of the
eigenvalue $\lambda$ on $\ep$, instead of Taylor expanding up to the
order 2, so that the group velocity becomes a pseudo-differential
operator $v(k+\ep\dx)$.

\paragraph{Rectification.} For general systems, even if the
initial data is purely oscillatory, a 
non-oscillating mode (or mean field) will be created. This is optical
\emph{rectification}. In \cite{JMR98}, Joly, M\'etivier and Rauch
allow such a phenomenon, still for non-dispersive systems. The case of
semilinear dispersive systems is treated by Lannes \cite{Lan98}. We
now illustrate the difficulty caused by interaction of oscillating and
non-oscillating modes, and simply look at the obtention of the first
profile equations. 

In the dispersive case, $\pi_0\neq \mbox{Id}$ (the origin belongs to
the characteristic variety $\mathcal{C}$ of $L$) and the set of nonzero 
eikonal frequencies, 
$$\mathcal{E} = \{ \alpha\in\Z^\star \mid
\alpha(-\omega,k)\in\mathcal{C} \},$$ 
is finite, in general. We assume that the only singular point of $C$
in $\{ \alpha(-\omega,k) \mid \alpha\in\mathcal{E} \}$ is the origin. 
Splitting profiles $U_n = \sum_{\alpha\in\Z} U_n^\alpha(T,X,t,x)~
e^{i\alpha\theta}$ into average $U_n^0$ and oscillations
$U_n^\star=U_n-U_n^0$, one gets from the first powers of $\ep$ in the
formal expansion of $L u^\ep -F(u^\ep)$,  
\begin{subequations} \label{fsteqprofdiffmoy}
\begin{align}
&\Pi ~ U_0^\star = U_0^\star , \qquad  \pi_0 ~ U_0^0 = U_0^0 ,
\label{polardiff} \\ 
& \pi_0 L_1(\partial_{t,x}) \pi_0 ~ U_0^0 = 0 , \label{transpmoydiff} \\
& \pi_0 L_1(\partial_{T,X}) \pi_0 ~ U_0^0 
- \pi_0 L_1(\partial_{t,x}) L_0^{-1} 
L_1(\partial_{t,x}) \pi_0 ~ U_0^0 \label{fsteqprofdiffmoyNLS} \\
& \qquad\qquad\qquad\qquad + \pi_0 F_0(U_0^\star,U_0^0) = 
-(\dt + \vec{v}_\alpha\cdot\dx) ~ U_1^\alpha , \notag 
\end{align}
\end{subequations}
and $\forall \alpha \in \Z \nz$, 
\begin{subequations} \label{fsteqprofdiffosc}
\begin{align}
& (\dt + \vec{v}_\alpha\cdot\dx) ~ U_0^\alpha = 0 ,
\label{transposcdiff} \\
& (\dT + \vec{v}_\alpha\cdot\dX) ~ U_0^\alpha 
- \pi_\alpha L_1(\partial_{t,x}) L(i\alpha\beta)^{-1} 
L_1(\partial_{t,x}) \pi_\alpha ~ U_0^\star \label{fsteqprofdiffoscNLS} \\
& \qquad\qquad\qquad\qquad + \pi_\alpha F_\alpha(U_0^\star,U_0^0) = 
-(\dt + \vec{v}_\alpha\cdot\dx) ~ U_1^\alpha. \notag 
\end{align}
\end{subequations}

As for the transport operator at the group velocity, an ``algebraic
lemma'' \cite{Tex04} shows that the second order operator 
$\pi_\alpha L_1(\partial_{t,x}) L(i\alpha\beta)^{-1}
L_1(\partial_{t,x}) \pi_\alpha$ is in fact $\frac{i}{2}
\partial_\xi^2\lambda(\alpha\beta)\cdot(\dx,\dx)$. On the other hand, 
the operator $\pi_0 L_1(\partial_{t,x}) \pi_0$, corresponding to low
frequencies, and called ``long-wave operator'', is symmetric
hyperbolic, but is not a transport operator, since the origin is a
singular point of the characteristic variety $\mathcal{C}$ of $L$. 
We simply know \cite{Lan98} that its characteristic variety 
$\mathcal{C}_{LW}$ is the tangent cone to $\mathcal{C}$ at the
origin. This cone may contain hyperplanes $\{ (\tau,\xi) \in \R^{1+d} 
\mid\tau+\vec{v}_\alpha\cdot\xi=0 \}$; this is precisely the 
rectification criterium. 

The difficulty now lies in the compatibility of the profile
equations. Lannes \cite{Lan98} defines nice average operators, which
are the analytic analogue of the algebraic projector $\Pi$. They
identify and separate the various propagation modes at scale $(T,X)$ 
and give necessary and sufficient
conditions for the profile equations to be solvable, with a
$t$-sublinear corrector $U_1$. Split the set of eikonal frequencies
$\mathcal{E}$ into the resonant set 
$$\mathcal{E}_R = \{ \alpha \in \mathcal{E} \mid \{
\tau+\vec{v}_\alpha\cdot\xi=0 \} \subset 
\mathcal{C}_{LW} \} , $$ 
--denoting $U_0^{0,\alpha}$ the corresponding
modes of $U_0^0$, and $U_0^{0,\alpha'}$, $\alpha'\in\mathcal{E}'$, the
others-- and the non-resonant set, $\mathcal{E}_{NR} =
\mathcal{E} \setminus \mathcal{E}_R$. Now, the abovementioned
conditions consist in replacing \eqref{fsteqprofdiffmoyNLS},
\eqref{fsteqprofdiffoscNLS}, for a resonant mode 
$\alpha\in\mathcal{E}_R$, with
\begin{subequations} \label{eqprofdiffRNLS}
\begin{align}
& \pi_0 L_1(\partial_{T,X}) \pi_0 ~ U_0^{0,\alpha} 
- \pi_0 L_1(\partial_{t,x}) L_0^{-1} 
L_1(\partial_{t,x}) \pi_0 ~ U_0^{0,\alpha} + \pi_0
F_0(U_0^\alpha,U_0^{0,\alpha}) = 0 , \label{eqprofdiffmoyRNLS} \\ 
& (\dT + \vec{v}_\alpha\cdot\dX) ~ U_0^\alpha - \frac{i}{2}
\partial_\xi^2(\alpha\beta)\cdot(\dx,\dx) ~ U_0^\alpha + \pi_\alpha
F_\alpha(U_0^\alpha,U_0^{0,\alpha}) = 0, 
\end{align}
\end{subequations}
whereas equations for non-resonant modes
($\alpha\in\mathcal{E}_{NR}, \alpha'\in\mathcal{E}'$) are decoupled
from the others,  
\begin{subequations} \label{eqprofdiffNRNLS}
\begin{align}
& \pi_0 L_1(\partial_{T,X}) \pi_0 ~ U_0^{0,\alpha'} 
- \pi_0 L_1(\partial_{t,x}) L_0^{-1} L_1(\partial_{t,x}) \pi_0
U_0^{0,\alpha'} + \pi_0 F_0(0,U_0^{0,\alpha'}) = 0 , \\
& (\dT + \vec{v}_\alpha\cdot\dX) ~ U_0^\alpha - \frac{i}{2}
\partial_\xi^2(\alpha\beta)\cdot(\dx,\dx) U_0^\alpha + \pi_\alpha
F_\alpha(U_0^\alpha,0) = 0. 
\end{align}
\end{subequations}
The corrector $\Pi U_1$ absorbs the difference between these 
equations and \eqref{fsteqprofdiffmoyNLS}, \eqref{fsteqprofdiffoscNLS},
\begin{equation*} 
\begin{split}
& \pi_0 L_1(\partial_{t,x}) \pi_0 ~ U_1^0 =  \pi_0 
\left[ F_0(U_0^\star,U_0^0) - \sum_{\alpha\in\mathcal{E}} 
F_0(U_0^\alpha,U_0^{0,\alpha}) \right], \\
& (\dt + \vec{v}_\alpha\cdot\dx) ~ U_1^\alpha = \pi_\alpha 
[ F_\alpha(U_0^\star,U_0^0) - F_\alpha(U_0^\alpha,U_0^{0,\alpha}) ], 
\quad \forall\alpha\in\mathcal{E}.
\end{split}
\end{equation*}

This mode decoupling produces a solvable system of profile equations
(polarization \eqref{polardiff}, linear hyperbolic evolution
\eqref{transpmoydiff}, \eqref{transposcdiff}, slow
nonlinear evolution \eqref{eqprofdiffRNLS}), \eqref{eqprofdiffNRNLS}), 
and clearly explains the rectification effect: 
even if $U_0^0$ vanishes initially, non-oscillating modes may be
created by nonlinear interactions in \eqref{eqprofdiffmoyRNLS}. The
asymptotics is then only at first order: the profile $U_0$ is
constructed, together with correctors $U_1$, $U_2$, whose secular
growth (see \cite{Lan03}) prevents from higher order asymptotics.

The analogue study for hyperbolic systems with variable coefficients
is done by Dumas \cite{Dum03a}, \cite{Dum04}, for WKB expansions with
several nonlinear phases: a rescaling in \eqref{singray} leads to a
``weakly nonplanar'' version of long-time diffraction,
\begin{equation*} 
u^\ep(t,x) \sim \ep^m \sum_{n\geq0}\ep^{n} U_n
(\ep t,\ep x,\vec{\psi}(\ep t,\ep x)/\ep,\vp(\ep t,\ep x)/\ep^2). 
\end{equation*}
The validation of such asymptotics necessitates coherence assumptions
on the phases $\vp$, $\vec{\psi}$ relatively to the operator $L$ as
well as to the tangent operators $\pi_\alpha L \pi_\alpha$.

\paragraph{Self-focusing.} Solutions to a (focusing)
nonlinear Schr\"odinger equation may blow up in finite time, and this
is usually interpreted in laser physics as the self-focusing
of the laser beam \cite{SS99}. Now, if the solution to the original
Mawxell equations describing the beam is globally defined, what does
this singularity mean? In \cite{Dum03b}, Dumas shows that the
diffractive profile blow-up corresponds to focusing for a perturbed
weakly nonlinear geometric optics model, and 
Schr\"odinger approximation is valid at least up to a time $t_\ep$ of the
order of a negative power of $\ln(1/\ep)$ before blow-up. This
shows that, even if each solution $u^\ep$ is defined globally in time,
it undergoes, between $t=0$ and $t=t_\ep$, an amplification by a
positive power of $\ln(1/\ep)$. 

\paragraph{Pulses and continuous spectra.} The propagation
of pulses over diffractive times is described by Alterman and Rauch 
in the linear \cite{AR00} and nonlinear case \cite{AR03}. The profiles
from \eqref{profdiff} then have compact support in $\theta$. An
important difference with diffractive optics for oscillatory
wavetrains resides in the profile equations
$$\dT \dz U - \frac{1}{2}\partial_\xi\lambda(\dz)(\dx^2+\dy^2)U = F(U).$$
The operator $\dz^{-1}$ is not defined on the space of profiles
considered, which can be seen, by Fourier transform, as a small
divisors problem. Alterman and Rauch solve this difficulty using
``infrared cut-offs'': they define approximate profiles by truncating
the low frequencies in the equation, and then show the convergence of
the approximate solution based on these approximate profiles. See also
\cite{SU03}: even if Schneider and Uecker do not compare their results
with the ones from the WKB method, they address the problem of existence
and stability of diffractive pulse solutions to nonlinear optics
Maxwell's equations, using a center manifold reduction (so that they
get exponential asymptotic stability, but only around some particular
family of solutions).

In \cite{BL02}, Barrailh and Lannes extend this approach to profiles
with continuous spectra, which model Raman scattering and lasers with
large spectrum.

\subsection{Transparency and larger amplitudes} \label{secstrongwaves}

\paragraph{Transparency.} The analysis above enlights the nonlinear
interactions leading to rectification. Unfortunately, computations on
physical models, such as Maxwell-Bloch equations , or ferromagnetic
Maxwell equations, 
reveal that the nonlinearities in \eqref{eqprofdiffmoyRNLS} simply
vanish! This phenomenon, called \emph{weak transparency}, is expressed as,
$$\forall U\in\C^N, \forall\alpha\in\Z , \quad \pi_\alpha \sum_{\beta\in\Z}
F_\alpha(\pi_\beta U,\pi_{\alpha-\beta} U) = 0.$$ 

In order to reach nonlinear regimes, one may then increase the
observation time, or the amplitude. The latter is studied (for
geometric optics time $\mathcal{O}(1)$) by Joly, M\'etivier, Rauch in 
\cite{JMR00b} for semilinear systems of Maxwell-Bloch type (including
the physical ferromagnetism system in space dimension one),
\begin{equation} \left\{
\begin{split}
& L(\ep\partial) u^\ep + \ep f(u^\ep,v^\ep) = 0 , \\
& M(\ep\partial) v^\ep + g(u^\ep,u^\ep) + \ep h(u^\ep,v^\ep) = 0 ,
\end{split} \right.
\end{equation}
with $L$ and $M$ symmetric hyperbolic operators as in \eqref{op}, and
$f$, $g$, $h$ bilinear nonlinearities. 
When weak transparency occurs, they look for solutions
with larger amplitude than the usual ones ($\mathcal{O}(1)$ instead 
of $\ep$). They show that weak transparency is
necessary to construct formal WKB expansions, and give a necessary and
sufficient criterium (strong transparency) for the stability of these
WKB solutions. Strong transparency and the particular structure of the
system in fact allow a nonlinear change of unknowns which brings back
to the weakly nonlinear setting. 

The same transparency property allows Jeanne \cite{Jea02} to construct
geometric optics asymptotics of large solutions to (semilinear)
Yang-Mills equations from general relativity. 

Using Joly-M\'etivier-Rauch's method, Colin \cite{Col02} rigorously 
derives Davey-Stewartson (DS) systems (coupling of (NLS) and a
hyperbolic or elliptic equation, modeling for example shallow-water, 
starting from Euler equations with free surface) from Maxwell-Bloch type
systems, over diffractive times. In space dimension one,
Schneider \cite{Sch98a}, \cite{Sch98b} obtains similar results \via
normal form techniques. Colin also shows that the obstruction for 
such a derivation from general hyperbolic systems comes from
rectification effects. Next, Colin and Lannes \cite{CL04} perform the
extension to general systems, and apply their results to
the Maxwell-Landau-Lifshitz ferromagnetism model (see Landau and Lifshitz
\cite{LL69}), getting mean-field generation. This 
corresponds to the physics papers of Leblond on pulse propagation,
deriving (NLS) \cite{Leb01}, combining with an expansion of optical
susceptibilities \cite{Leb02}, or deriving (DS) in ferromagnetic media
\cite{Leb96}, \cite{Leb99}. 

For ill-prepared data allowing rectification, they build a 
(DS) approximation (valid only over times of the order of 
$\ln(1/\ep)$). To this end, they need a long-wave correction to the 
original Ansatz. Furthermore, they assume (in order for the limit
(DS) system to be well-posed, \ie to be a
Schr\"odinger-elliptic coupling) that long-wave--short-wave
resonance (CROLOC, in French) does not occur, which means that the
characteristic varieties of the long-wave operator and of the original
operator are not tangent away from the origin. In order to reach
$\mathcal{O}(1/\ep)$ times with rectification, they consider
\cite{CL01} solutions with size $\rep$ (intermediate between $1$ and
$\ep$), and use four-scale profiles to obtain their CROLOC. A
numerical study of long-wave--short-wave resonance is due to Besse and
Lannes \cite{BL01}.

Some works on waves in plasmas rely on similar techniques. The basic
model is then the quasilinear Euler-Maxwell system (see Sulem and
Sulem \cite{SS99}). Colin, Ebrard, Gallice and Texier \cite{CEGT04}
study a simpler model, a Klein-Gordon-wave coupling, for which a
change of unknowns leads to a semilinear system with weak
transparency property. The diffractive time approximation is the
Zakharov system, for electromagnetic field $u$, and ion population
$n$, 
\begin{equation*} \left\{
\begin{split}
& i \dt u + \Delta u = n u ,\\
& \dt^2 n - \Delta n = \Delta |u|^2. 
\end{split} \right.
\end{equation*}
For the full Euler-Maxwell system, Texier \cite{Tex05}
rigorously derives, in the geometric optics regime, a weak form of the
above system (where $\ep$ stands in front of the terms $\Delta u$ and
$nu$, thanks to ``generalized WKB asymptotics'', \ie without
eliminating the residual from profile equations.

\paragraph{Conservation laws with a linearly degenerate field.}
Cheverry, Gu\`es and M\'etivier have adapted the notion of
transparency above to quasilinear systems of conservation laws. They
classify nonlinear regimes on times $\mathcal{O}(1)$ with asymptotics
$$u^\ep(t,x) \sim u_0(t,x) + \sum_{k\geq1} \ep^{k/l} U_k(t,x,\ph/\ep).$$
The value $l=1$ corresponds to weakly nonlinear regime; $l=2$, to
\emph{strong oscillations}; $l\geq3$, to \emph{turbulent
  oscillations}; $l=\infty$, to \emph{large amplitude} (\ie
$\mathcal{O}(1)$) oscillations (more singular solutions are the
stratified solutions of Rauch and Reed \cite{RR88}, studied in the
quasilinear case by Corli and Gu\`es \cite{CG01}). A difficulty in WKB
analysis, when dealing with amplitudes larger than the weakly nonlinear
ones, comes from the fact that the hierarchy of profile equations
changes. In particular, a coupling appears between phase and amplitude
in the eikonal equation. 

In \cite{CGM03}, they show that linear degeneracy of one of the fields is a
weak transparency condition that ensures existence of formal strong
solutions. In \cite{CGM04}, they study existence and stability of
large amplitude (polarized) waves, under a stronger transparency
assumption (existence of a good symmetrizer and constant multiplicity
of the linearly degenerate eigenvalue) that relates the system to
Euler non-isentropic gas dynamics equations --see M\'etivier and
Schochet \cite{MS01} about the stability of large oscillations in time
for the entropy. Earlier results on this topic were only formal ones
(Serre \cite{Ser95}), or in space dimension one: see Peng \cite{Pen92},
Heibig \cite{Hei93}, Corli and Gu\`es \cite{CG01}, and Museux \cite{Mus04}.

Cheverry continues this work, in connection with turbulence. 
In \cite{Che04a}, he proves that (vanishing) viscosity can
compensate the lack of transparency. In \cite{Che04b}, he converts the
stability problem into a ``cascade of phases'' phenomenon, where the
phase has an asymptotic expansion, whose coefficients are coupled with
amplitudes.

\paragraph{Supercritical WKB solutions to (NLS).} 
Coupling between phase and amplitude also occurs for large amplitude
oscillating solutions to nonlinear Schr\"odinger equations. Existence 
of such solutions before caustics is obtained by Grenier \cite{Gre98} 
for a single nonlinear phase $\ph$, with infinitely accurate WKB 
expansions (extending results of G\'erard \cite{Ger93}), thanks to a 
``perturbed phase'' technique, $\ph \sim \sum_k \ep^k \ph_k$.

\section{Long waves} \label{seclong}

We finally mention the situation \emph{a priori} opposite to highly
oscillatory problems, where wavelength is large. This is the typical
framework of water waves, modelled by Euler equations with free
surface, (EFS); see \cite{SS99}. However, here again come into play the
coupling with a mean field, and the ``long wave operator'', as for
diffractive optics rectification effects from previous section. This
explains, for example, the formation of pairs of waves travelling in
opposite directions. 

Due to symmetries, the second-order differential terms vanish, in
long-time asymptotic models, and the limit equations take the form of
Korteweg-de Vries (KdV) equations. For one-dimensional quadratic
quasilinear dispersive systems of type \eqref{syst}, the appropriate
Ansatz is
$$u^\ep(t,x) \sim \ep^2 \sum_{n\geq0} \ep^n U_n(\ep^2t,t,x),$$
and profiles satisfy transport equation at scale $(t,x)$ so that,
$$U_0(\ep^2t,t,x) = V_0(T,y)_{|{T=\ep^2 t, y=x\pm t}}.$$ 
Furthermore, they are solutions to (KdV) in variables $(T,y)$,
$$\dT V_0 + \frac{1}{6} \dy^3 V_0 + \frac{3}{4} \dy(V_0^2) = 0. $$
Ben Youssef and Colin \cite{BYC00}, as well as Schneider
\cite{Sch98c}, derive this asymptotics for simplified general systems,
getting decoupled equations for each mode. The case of (EFS) is
treated by Schneider and Wayne with \cite{SW02} or without \cite{SW00} 
surface tension. 
The case of transverse perturbations is more singular, since it leads
(for profiles $U(\ep^2 t,\ep x_2,t,x_1)$) to a
Kadomtsev-Petviashvili-type equation (KP),
$$\dT V + \frac{1}{6} \partial_y^3 V + \frac{1}{2} \partial_y^{-1}
\dY^2 V + \frac{3}{4} \dy(V^2) = 0. $$
This is obtained by Gallay and Schneider for unidirectional waves, and
then by Ben Youssef and Lannes \cite{BYL02} for pairs of waves. Note
that (KP) equations are singular with respect to low frequency (in $y$)
solutions. The ``infrared cut-off'' technique from section \ref{seclongtime}
is thus useful here --and this singularity poses the problem of
consistency of the approximate solution, since residuals
$Lu^\ep-F(u^\ep)$ may not be small (see Lannes \cite{Lan03b}).

The equations are asymptotically decoupled, but coupling appears at
least when propagation takes place in a bounded domain (or for
periodic solutions). Taking coupling into account also improves
convergence rates. Ben Youssef and Lannes highlight such coupling effects
between (KP) equations simply using ``generalized WKB expansions'' as
described in paragraph \ref{seclongtime}. A well-known coupled 
approximation of (EFS) is the Boussinesq system. In fact, Bona, Colin
and Lannes \cite{BCL03} obtain, in 2 or 3 space dimension, a result
ensuring the simultaneous validity (or non-validity) of a whole
three-parameters family of such systems, including Boussinesq. 
For numerics on this topic, see Labb\'e and Paumond \cite{LP04}.

The same regime is of interest in micromagnetism, \ie for
Maxwell-Landau-Lifshitz system (where the medium responds to the
electromagnetic field \via magnetization). Colin, Galusinski and Kaper
\cite{CGK02} study the propagation of pairs of travelling waves in space
dimension one, and derive a semilinear heat equation, whereas for
two-dimensional perturbations, Sanchez \cite{San05} obtains 
Khokhlov-Zabolotskaya equations.


\bibliographystyle{alpha}
\bibliography{reviewGO}

\end{document}